\documentclass[12pt]{amsart}
\usepackage{amsmath}
\usepackage{amsfonts}
\numberwithin{equation}{section}

\newtheorem{theorem}{Theorem}[section]
\theoremstyle{plain}

\newtheorem{corollary}[theorem]{Corollary}

\newtheorem{lemma}[theorem]{Lemma}

\newtheorem{proposition}[theorem]{Proposition}

\theoremstyle{definition}

\theoremstyle{remark}
\newtheorem{remark}[theorem]{Remark}

\input{tcilatex}

\def\C{{\mathbb C}}

\def\R{{\mathbb R}}
\def\B{{\mathbb B}}
\def\H{{\mathbb H}}

\def\bn{{\B^n}}

\def\de{\partial}

\def\id{{\sf id}}

\def\Aut{{\sf Aut}}

\def\Re{\func{Re}}

\def\v{\varphi}

\def\D{\mathbb D}

\def\la{\langle}
\def\ra{\rangle}

\def\lfmb{{\rm LFM}(\bn,\bn)}

\begin{document}
\title[Infinitesimal generators]{Infinitesimal generators associated with
semigroups of linear fractional maps}
\author[F. Bracci]{Filippo Bracci}
\address{F. Bracci: Dipartimento Di Matematica\\
Universit\`{a} di Roma \textquotedblleft Tor Vergata\textquotedblright\ \\
Via Della Ricerca Scientifica 1, 00133 \\
Roma, Italy} \email{fbracci@mat.uniroma2.it}
\thanks{}
\author[M.D. Contreras]{Manuel D. Contreras$^\dag$}
\address{M.D. Contreras \and  S. D\'{\i}az-Madrigal: Camino de los Descubrimientos, s/n\\
Departamento de Matem\'{a}tica Aplicada II \\
Escuela Superior de Ingenieros\\
Universidad de Sevilla\\
41092, Sevilla\\
Spain.} \email{contreras@esi.us.es, madrigal@us.es}
\author[S. D\'{\i}az-Madrigal]{Santiago D\'{\i}az-Madrigal$^\dag$}

\thanks{$^\dag$Partially supported by the \textit{Ministerio
de Ciencia y Tecnolog\'{\i}a} and the European Union (FEDER) project
BFM2003-07294-C02-02 and by \textit{La Consejer\'{\i}a de
Educaci\'{o}n y Ciencia de la Junta de Andaluc\'{\i}a.}}

\begin{abstract} We characterize the infinitesimal generator of a
semigroup of linear fractional self-maps of the unit ball in
$\mathbb C^n$, $n\geq 1$. For the case $n=1$ we also completely
describe the associated Koenigs function and we solve the embedding
problem from a dynamical point of view, proving, among other things,
that a generic semigroup of holomorphic self-maps of the unit disc
is a semigroup of linear fractional maps if and only if it contains
a linear fractional map for some positive time.
\end{abstract}

\maketitle

\section*{Introduction}

Let $U$ be an open set in $\C^n$, $n\geq 1$. A continuous semigroup
$(\varphi_t)$ of holomorphic functions in $U$ is a continuous
homomorphism from the additive semigroup of non-negative real
numbers into the composition semigroup of all holomorphic self-maps
of $U$  endowed with the compact-open topology. In other words, the
map $[0,+\infty )\ni t\mapsto (\varphi_{t})\in\mathrm{Hol}(U,U)$
satisfies the following three conditions:

\begin{enumerate}
\item $\varphi _{0}$ is the identity map $\id_U$ in $U,$

\item $\varphi _{t+s}=\varphi _{t}\circ \varphi _{s},$ for all
$t,s\geq 0,$

\item $\varphi _{t}$ tends to $\id_U$ as $t$ tends to $0$
uniformly on compacta of $U$.
\end{enumerate}

Any element of the above family $(\varphi _{t})$ is called an {\sl
iterate} of the semigroup. If the morphism $\lbrack 0,+\infty )\ni
t\mapsto \varphi _{t}\in \mathrm{Hol}(U, U)$ can be extended
linearly to $\mathbb{R}$ (and then necessarily each $\varphi_t$ is
invertible) we have a group  of automorphisms of $U$. It is well
known (see, {\sl e.g.}, \cite{Abate} or \cite{Shoiket}) that every
iterate of a semigroup is injective and that if for some $t_0>0$
the iterate $\v_{t_0}\in \Aut(U)$ then $\v_{t}\in \Aut(U)$ for all
$t>0$ and the semigroup is indeed a group of automorphisms of $U$.

The theory of semigroups in the unit disc $\mathbb{D}$ has been
deeply studied and applied in many different contexts. We refer
the reader to the excellent monograph by Shoikhet \cite{Shoiket}
and references therein for more information about this.

For several and natural reasons, those semigroups in $\mathbb{D}$
whose iterates belong to some concrete class of holomorphic
functions relatively easy to handle are of special interest.
Undoubtedly, the most important example of this kind is that of
semigroups of linear fractional self-maps of $\D$ (shortly,
$\mathrm{LFM}(\mathbb{D},\mathbb{D})$), namely, those semigroups
$(\v_t)$ in $\mathbb{D}$ for which   every iterate $\varphi _{t}\in
\mathrm{LFM}(\mathbb{D},\mathbb{D})$. Among those semigroups of
linear fractional maps, there are  the groups of automorphisms of
$\D$.

Semigroups in $\mathbb{D}$ can be first classified looking at
their common Denjoy-Wolff fixed point. It can be proved  (see
\cite{Siskakis-tesis}, \cite{Siskakis-JLMS},
\cite{Contreras-Diaz-Madrigal-analytic-flows}) that any semigroup
$(\varphi _{t})$ in $\mathbb{D}$ belongs to one and only one of
the following five classes :

\begin{enumerate}
\item   {\sl trivial-elliptic}:  all the iterates are the identity
$\id_{\mathbb{D}}$ map. \item {\sl Neutral-elliptic}: there exists
$\tau \in \mathbb{D}$ with $\varphi _{t}(\tau )=\tau $ and
$\left\vert \varphi _{t}^{\prime }(\tau )\right\vert =1$ for every
$t>0$ and $\v_{t}\neq \id_\D$ for some $t>0$. \item {\sl
Attractive-elliptic}: there exists $\tau \in \mathbb{D}$ with
$\varphi _{t}(\tau )=\tau $ and $\left\vert \varphi _{t}^{\prime
}(\tau )\right\vert <1,$ for every $t>0.$ \item {\sl Hyperbolic}:
there exists $\tau \in \partial \mathbb{D}$ with $\lim_{r\to
1^-}\varphi _{t}(r\tau )=\tau $ and $\lim_{r\to 1^-}\varphi
_{t}^{\prime }(r\tau )<1,$ for every $t>0.$ \item {\sl Parabolic}:
there exists $\tau \in
\partial \mathbb{D}$ with $\lim_{r\to 1^-}\varphi _{t}(r\tau )=\tau $ and $\lim_{r\to 1^-}\varphi
_{t}^{\prime }(r\tau )=1,$ for every $t>0$.
\end{enumerate}

In cases $(2)$ to $(5)$, the   point $\tau $ is unique and it is
called the {\sl Denjoy-Wolff point} of the semigroup.  According
to the Julia-Wolff-Carath\'eodory theorem (see, {\sl e.g.}
\cite{Abate}, \cite{Shapiro-libro}), in cases (4) and (5) all
iterates have non-tangential (or angular) limit $\tau$ at $\tau$
and their first derivatives have a non-tangential limit at $\tau$
given by a real number in $(0,1]$. As customary, we call {\sl
elliptic semigroup} any semigroup in the classes (1) to (3) and
{\sl non-elliptic} the semigroups in the classes (4) and (5). An
analogous classification is available for semigroups of $\B^n$,
the unit ball of $\C^n$ (see, {\sl e.g.}, \cite{Abate} or
\cite{Bracci-Con-Diaz} or Section two).

To any semigroup  in $\mathbb{D}$  there are attached two analytic
objects which can be used to describe the dynamical behavior: the
infinitesimal generator and the Koenigs map. We are going to quickly
recall how they are defined.

Given a semigroup $(\varphi _{t})$ in $\mathbb{D}$, it can be proved
(see \cite{Shoiket}, \cite{Berkson-Porta}) that there exists a
unique holomorphic function $G:\mathbb{D\rightarrow C}$ such that,
for each $z\in \mathbb{D},$ the trajectory
\begin{equation*}
\gamma _{z}:[0,+\infty )\rightarrow \D,\quad \quad t\mapsto \gamma
_{z}(t):=\varphi _{t}(z)
\end{equation*}%
is the solution of the Cauchy problem
\begin{equation*}
\left\{
\begin{array}{l}
\dot{w}=G(w) \\
w(0)=z. %
\end{array}%
\right.
\end{equation*}%
Moreover, $G(z)=\lim_{t\rightarrow 0^{+}}(\varphi _{t}(z)-z)/t$, for
every $z\in \mathbb{D}.$ The function $G$ is called the {\sl
infinitesimal generator} (or the {\sl semi-complete vector field})
of $(\varphi _{t}).$
There is a very nice representation, due to Berkson and Porta \cite%
{Berkson-Porta}, of those holomorphic functions of the disc which
are infinitesimal generators. Namely:

\begin{theorem}[Berkson-Porta]\label{24} A  holomorphic function $G:\mathbb{D\rightarrow C}$ is the
infinitesimal generator of a semigroup in $\mathbb{D}$ if and only
if there exists a point $b\in \overline{\mathbb{D}}$ and a
holomorphic function $p:\mathbb{D\rightarrow C}$ with
$\func{Re}p\geq 0$ such that
\begin{equation*}
G(z)=(z-b)(\overline{b}z-1)p(z),\text{ \quad }z\in \mathbb{D}.
\end{equation*}
\end{theorem}

The point $b$ in Berkson-Porta's theorem is exactly the Denjoy-Wolff
point of the semigroup, unless the semigroup is  trivial. Other
alternative descriptions of infinitesimal generators can be found in
\cite[Section 3.6]{Shoiket} (where a different sign convention is
chosen).

As for the Koenigs function,  if $(\varphi _{t})$ is a semigroup in
$\mathbb{D}$ with Denjoy-Wolff point $\tau \in \mathbb{D}$ and
infinitesimal generator $G,$ then there exists a unique
univalent function $h\in \mathrm{Hol}(\mathbb{D},\mathbb{C})$ such that $%
h(\tau )=0,$ $h^{\prime }(\tau )=1$ and, for every $t\geq 0,$%
\begin{equation*}
h\circ \varphi _{t}(z)=\frac{d\varphi_{t}}{d z}(\tau
)h(z)=e^{G^{\prime }(\tau )t}h(z).
\end{equation*}%
While, if $(\varphi _{t})$ is a semigroup in $\mathbb{D}$ with
Denjoy-Wolff point $\tau \in \partial \mathbb{D},$ then there exists
a unique univalent function $h\in
\mathrm{Hol}(\mathbb{D},\mathbb{C})$ such
that $h(0)=0$ and, for every $t\geq 0,$%
\begin{equation*}
h\circ \varphi _{t}(z)=h(z)+t.
\end{equation*}%
In both cases, the function $h$ is called the {\sl Koenigs function}
of the semigroup $(\varphi _{t})$.

Infinitesimal generators can be defined also in several variables
(see \cite{Abate}, \cite{AERS}, and \cite{Shoiket}) even if their
characterizations are not so easy to handle as in the one
dimensional case. However the construction of Koenigs' functions
in several variables is still at a pioneeristic step and it has
been successfully developed only for linear fractional semigroups
(see \cite{Bracci-Con-Diaz} and \cite{Cowen-Others} for a
different construction but only for a single linear fractional map
of $\B^n$) and, in the realm of discrete iteration, for regular
hyperbolic self-maps of $\B^n$ (see \cite{Bracci-Gentili}).

A related problem is that of embedding a given holomorphic self-map
of a domain $U\in\C^n$ into a semigroup of holomorphic self-maps of
$U$. Such a problem has been studied since long, see \cite{Shoiket}
for a good account of available results.

Semigroups of $\Aut(\mathbb{D})$ are quite well-understood (see
\cite {Begson-Porta-TAMS}) in terms of infinitesimal generators
and Koenigs functions.  Moreover,  the embedding problem is
completely solved: every automorphism of the unit disc can be
embedded in a suitable (semi)group of $\Aut(\mathbb{D})$.

However, strange as it may seem, the situation for semigroups of
$\mathrm{LFM}(\B^n,\B^n)$ (even for $n=1$) is still not completely
clear. Despite the fact that   convergence questions are basically
understood, there is not a full description of their basic
theoretical elements such as infinitesimal generators and Koenigs
functions. The corresponding embedding problem for $n=1$ has been
treated in the literature and it is known that, in general, the
answer is negative. In most cases   some analytic criteria for
deciding the solvability of this problem  are available (see
\cite[Sections 4.3 and 5.9]{Shoiket}) but the problem in its
complete generality seems to be still open.

In this paper, we consider those problems for semigroups of
$\mathrm{LFM}(\B^n,\B^n)$,  especially for $n=1$. We completely
characterize infinitesimal generators of semigroups of linear
fractional self-maps of the ball  (in one and several variables). In
particular we prove the following result:

\begin{theorem}\label{main} Let $(\v_t)$ be a semigroup of holomorphic self-maps of
$\B^n$. Then $(\v_t)$ is a semigroup of linear fractional maps if
and only if there exist $a,b\in \C^n$ and $A\in \C^{n\times n}$ (not
all of them zero) such that the infinitesimal generator of $(\v_t)$
is
\[
G(z)=a-\la z,a\ra z-[Az+\la z,b\ra z]
\]
with $|\la b,u\ra | \leq \Re \la Au, u\ra$, for all $u\in \partial
\B^n$.
\end{theorem}

For the case $n=1$, in Section two we present  a more precise
statement classifying all possible cases (see Propositions
\ref{poldeg1} and \ref{poldeg2}), and we present a precise
description of Koening maps for linear fractional semigroups (see
Proposition \ref{kperlfm}). Finally, in the third section we deal
with the embedding problem. First we prove the following
rigidity result:
\begin{theorem}
Let $(\v_t)$ be a semigroup in $\D$. If for some $t_0>0$ the iterate
$\v_{t_0}$ is a linear fractional self-map of $\D$ then $\v_t$ is a
linear fractional self-map of $\D$ for all $t\geq 0$.
\end{theorem}
Then we settle the embedding problem for a linear fractional
self-map  of $\D$ proving that it can be embedded in a semigroup of
$\D$ if and only if it can be embedded into a semigroup of linear
fractional self-maps of $\D$ if and only if it can be embedded into
a group of M\"obius transformations of the Riemann sphere
$\C_\infty$ which for $t\geq 0$ preserves the unit disc (see Theorem
\ref{embedding}). Finally, we give a simple criterium for
embeddability of linear fractional self map of $\D$ (see Proposition
\ref{embedd}).

\bigskip

Part of this research has been carried out while the second and the
third quoted authors were visiting the University of Florence. These
authors want to thank the Dipartimento di Matematica
\textquotedblleft U. Dini\textquotedblright, and especially
professor G. Gentili, for   hospitality and support.

\section{Infinitesimal Generators in Several Variables}

Following \cite{Cowen-McClauer}, we say that a map $\varphi
:\mathbb{B}^{n} \to \C^n$ is a {\sl linear fractional map} if there
exist a complex $n\times n$ matrix $A\in \C^{n\times n}$, two column
vectors $B$ and $C$ in $\mathbb{C}^{n},$ and a complex number $D\in
\mathbb{C}$ satisfying
\begin{equation*}
{\rm (i)}\text{ }\left\vert D\right\vert >\left\Vert
 C\right\Vert ;\quad \text{ }{\rm (ii)}%
\text{ }DA\neq BC^{\ast },
\end{equation*}%
such that
\begin{equation*}
\varphi (z)=\frac{A z+B}{\left\langle z,C\right\rangle +D},\text{ \
}\quad z\in \mathbb{B}^{n}.
\end{equation*}%
Condition ${\rm (i)}$ implies that $\left\langle z,C\right\rangle
+D\neq 0$ for every $z\in \mathbb{B}^{n}$ and therefore, $\varphi $
is actually
holomorphic in a neighborhood of the closed ball. In fact, $\varphi \in \mathrm{Hol}(r%
\mathbb{B}^{n};\mathbb{C}^{n})$ for some $r>1.$ On the other hand,
condition ${\rm (ii)}$ just says that $\varphi $ is not constant. If
the image $\varphi(\bn)\subset\bn$, then we say that $\varphi$ is a
linear fractional self-map of $\bn$ and write $\varphi\in\lfmb$.

It is worth recalling that if $\varphi \in \lfmb$ has no fixed points in $%
\mathbb{B}^{n}$, then there exists a unique point $\tau \in
\partial \mathbb{B}^{n}$ such that $\varphi (\tau )=\tau $ and
$\la d\varphi _{\tau }(\tau ),\tau \ra=\alpha(\v)$ with
$0<\alpha(\v)\leq 1$ (see, {\sl e.g.}, \cite[Theorem
2.2]{Bisi-Bracci}). We call $\tau$ the {\sl Denjoy-Wolff point} of
$\v$ and $\alpha(\v)$ the {\sl boundary dilatation coefficient} of
$\v$.

A semigroup in $\mathrm{Hol}(\bn ;\bn )$ is a {\sl semigroup of
linear fractional maps} if  $\v_t\in \lfmb$ for all $t\geq 0$.

Likewise the one-dimensional case, given a semigroup $(\v_t)$ in
$\mathrm{Hol}(\B^n ,\B^n )$ there exists a holomorphic map
$G\in\mathrm{Hol}(\B^n ;\C^n )$, the {\sl infinitesimal generator}
of the semigroup, such that
$$
\frac{\partial \v_t}{\partial t}=G\circ \v_t
$$
for all $t\geq 0$ (see, {\sl e.g.}, \cite[Proposition
2.5.22]{Abate}).

The following result, essentially due to Abate \cite{Abate} (see
also \cite{AS} and \cite{Bracci-Con-Diaz}) allows to talk about
elliptic, hyperbolic, and parabolic semigroups in $\B^n$:

\begin{theorem}
Let $(\varphi _{t})$ be a semigroup in $\mathrm{Hol}%
(\mathbb{B}^{n},\mathbb{B}^{n}).$ Then, either all the iterates have
a common fixed point in $\mathbb{B}^{n}$ or all the iterates
$\varphi _{t}$ $(t>0)$ have no fixed points in $\mathbb{B}^{n}$ and
then they share the same Denjoy-Wolff point $\tau\in \de \bn$.
In this case, there exists $0<r\leq 1$ such that $\alpha _{t}=r^{t}$, where $%
\alpha _{t}:=\alpha(\v_t)$ denotes the boundary dilatation
coefficient of $\varphi _{t}$ (for $t>0$) at $\tau$.
\end{theorem}

A detail study of semigroups of linear fractional maps in $\mathrm{Hol}%
(\mathbb{B}^{n},\mathbb{B}^{n})$ can be found in
\cite{Bracci-Con-Diaz}. In this section, we present a
characterization of holomorphic functions $G:\B^n \to \C^n$ which
are infinitesimal generators of semigroups of linear fractional
maps. Before that, we need the following technical lemma.

\begin{lemma}\label{JWC-for-infin}
Let $(\v_t)$ be a semigroup of holomorphic self-maps of $\B^n$ with
associated infinitesimal generator $G$. Assume that $e_1\in \de
\B^n$ is the Denjoy-Wolff point of $(\v_t)$ and that $G$ extend
$C^1$ in a open neighborhood of $e_1$. Then
\begin{enumerate}
\item $G(e_1)=0$.
\item $\la dG_{e_1}(e_1), e_1\ra \in \R$.
\item $\la dG_{e_1}(e_j), e_1\ra =0$ for all $j=2,\ldots, n$.
\end{enumerate}
\end{lemma}
\begin{proof}
>From $G(\v_t(z))=\frac{\de}{\de t}\v_t(z)$ it follows that $\v_t$
are $C^1$ in a neighborhood of $e_1$ as well, and  then (1) follows.
Also for any $v\in \C^n$, we have $\la
dG_{e_1}(d(\v_t)_{e_1})v,e_1\ra=\frac{\de}{\de t}\la d(\v_t)_{e_1}v,
e_1\ra$ and then (2) and (3) follow from Rudin's version of the
classical Julia-Wolff-Carath\'eodory theorem in $\B^n$ applied to
$\v_t$, see \cite{RU} or \cite{Abate}.
\end{proof}

In the proof of Theorem \ref{main} we will  also make use of the
following generalization of Berkson-Porta's criterion due to
Aharonov, Elin, Reich and Shoikhet (see \cite[Theorem 4.1]{AERS},
where however a different sign convention is chosen because they
look at the problem $\frac{\partial \v_t}{\partial t}=-G\circ
\v_t$):

\begin{theorem}\label{shiokhet}
Let $F\in\mathrm{Hol}(\B^n,\C^n)$. Then $F$ is the infinitesimal
generator of a semigroup of holomorphic self-maps of   $\B^n$ fixing
the origin   if and only if $F(z)=-Q(z)  z$ where $Q(z)$ is a
$(n\times n)$-matrix with holomorphic entries such that
\[
\Re \langle Q(z)z,z\rangle\geq 0,
\]
for all $z\in\B^n$.
\end{theorem}

And now we can prove Theorem \ref{main}:

\begin{theorem} Let $(\v_t)$ be a semigroup of holomorphic self-maps of
$\B^n$. Then $(\v_t)$ is a semigroup of linear fractional maps if
and only if there exist $a,b\in \C^n$ and $A\in \C^{n\times n}$ (not
all of them zero) such that the infinitesimal generator of $(\v_t)$
is
\[
G(z)=a-\la z,a\ra z-[Az+\la z,b\ra z]
\]
 with
\begin{equation}\label{constrain}
|\la b,u\ra | \leq \Re \la Au, u\ra ,
\end{equation}
for all $u\in \partial \B^n$.
\end{theorem}
\begin{proof}
Suppose first that $(\v_t)$ is a semigroup of linear fractional maps
given by
\[
\v_t(z)=\frac{A_tz+B_t}{\la z, C_t\ra +1},
\]
for some holomorphic functions $t\mapsto A_t\in \C^{n\times n}$,
$t\mapsto B_t, C_t\in \C^n$. Differentiating with respect to $t$ at
$t=0$ and denoting by $A=\frac{\de A_t}{\de t}|_{t=0}$, $B=\frac{\de
B_t}{\de t}|_{t=0}$, $C=\frac{\de C_t}{\de t}|_{t=0}$ we obtain
\[
G(z):=\frac{\de \v_t}{\de t}|_{t=0}=B-\la z, B\ra z-[Az+\la z,
B+C\ra z].
\]
Thus the infinitesimal generator $G$ of $(\v_t)$ has the required
expression. We only need to verify \eqref{constrain}. To this aim,
 if $B\neq 0$ then $H(z)=-B+\la
z,B\ra z$ is the infinitesimal generator of a  group of  hyperbolic
automorphisms of $\B^n$ (this follows either by a direct simple
computation or by applying \cite[Theorem 3.1]{AERS}). Since the set
of infinitesimal generators is a (real) cone (see, for example,
\cite[Corollary~2.5.59]{Abate}) it follows that $G+H$ is an
infinitesimal generators of a semigroup of holomorphic self-maps of
$\B^n$. Let $b=B+C$. By Theorem~\ref{shiokhet} applied to $G+H$, we
have
\begin{equation*}
\Re \la Az, z\ra + \Vert z\Vert ^2 \Re \la z,b\ra  \geq 0 ,
\end{equation*}
for all $z\in  \B^n$. Now, writing $z=ru$ for $u\in \partial   \B^n$
and   $r\in [0,1)$, it follows that
\begin{equation*}
\Re \la Au, u\ra + r \Re \la u,b\ra  \geq 0 \quad \forall
u\in\de\B^n,\ r\in [0,1),
\end{equation*}
which implies \eqref{constrain}.

Conversely, suppose $G(z)=a-\la z,a\ra z-[Az+\la z,b\ra z]$
satisfies \eqref{constrain}. Then, we have that
\begin{equation}\label{constrain2}
\Re \la Az, z\ra + \Vert z\Vert ^2 \Re \la z,b\ra  \geq 0 ,
\end{equation}
for all $z\in  \B^n$. We can write $G(z)=H(z)+P(z)$ with $H(z)=a-\la
z,a\ra z$ and $P(z)=-[Az+\la z,b\ra z]$. As before, one can prove
that $H$ is the infinitesimal generator of a (semi)group of
holomorphic maps in $\B^n$. By \eqref{constrain} and Theorem
\ref{shiokhet} the function $P(z)$ is an infinitesimal generator of
a semigroup of holomorphic self-maps of $\B^n$ as well. Therefore
$G$ is an infinitesimal generator of a semigroup $(\v_t)$ of
holomorphic self-maps of $\B^n$. The point now is to show that such
a semigroup is composed by linear fractional maps.  There are two
cases: either $G(z_0)=O$ for some $z_0\in \B^n$ or $G(z)\neq O$ for
all $z\in \B^n$.

If $G(z_0)=O$ then $\v_t(z_0)=z_0$ for all $t\geq 0$. In this case
we first rotate conjugating with a unitary matrix $U$ in order to
map $z_0$ to the point $r e_1$ where $r:=\Vert z_0\Vert <1$. In
terms of infinitesimal generators this amounts to send $G$ to
$UGU^\ast$, and \eqref{constrain2} is preserved. Next, we move $r
e_1$ to the origin $0$ by means of the transform $T\in \Aut(\B^n)$
\[
T(\zeta,w)=\frac{(r-\zeta)e_1-(1-r^2)^{1/2}(0,w)}{1-r\zeta}, \quad
(\zeta,w)\in \C\times \C^{n-1}.
\]
The automorphism $T$ is an involution, that is, $T(T(z))=z$ for all
$z\in \B^n$ (see \cite[Section 2.2.1]{Abate} for this and others
properties of this $T$). The infinitesimal generator $G$ is sent to
$\tilde{G}(z)=dT_{T(z)} G(T(z))$. Let $\delta:=(1-r^2)^{1/2}$. A
direct computation shows that (with obvious notations)
\begin{equation*}
\begin{split}
dT_{T(\zeta,w)} G(T(\zeta,w))&=\delta^{-2}\left(
                               \begin{array}{cc}
                                 1-r\zeta & 0 \\
                                 -rw & \delta {\sf Id} \\
                               \end{array}
                             \right)A\left(
                                       \begin{array}{c}
                                         r-\zeta \\
                                         -\delta w \\
                                       \end{array}
                                     \right)\\&+\frac{\la (r-\zeta, -\delta w), b\ra}{1-r\zeta}\left(
                                                                                             \begin{array}{cc}
                                                                                               1-r\zeta & 0 \\
                                                                                               -rw & -\delta {\sf Id} \\
                                                                                             \end{array}
                                                                                           \right)\left(
                                       \begin{array}{c}
                                         r-\zeta \\
                                         -\delta w \\
                                       \end{array}
                                     \right).
\end{split}
\end{equation*}
Then $\tilde{G}$ is of the form $\tilde{G}(z)=-[Mz+\la z, c\ra z]$
for some $M\in \C^{n\times n}$ and $c\in \C^n$ and, since it is an
infinitesimal generator, again  by Theorem \ref{shiokhet}, it
satisfies
\begin{equation}\label{c2}
\Re \la Mz, z\ra +\|z\|^2\Re \la z,c\ra \geq 0
\end{equation}
for all $z\in \B^n$. Notice that $-M$ is dissipative because,
looking at degrees in $z$, equation~\eqref{c2} implies that $\Re \la
Mz, z\ra\geq 0$ for all $z\in \B^n$.

Assume first that $-M$ is not asymptotically stable. Thus there
exists $v\in \C^n$, $\|v\|=1$, such that $Mv=i\alpha v$ for some
$\alpha\in \R$. Substituting $z=\zeta v$ for $\zeta\in\C$,
$|\zeta|<1$ in \eqref{c2} we obtain that $\la v, c\ra=0$. Now let
$z=\zeta v + \epsilon c$ for $|\zeta|<1$ and $\epsilon\in\R$ small.
Substituting this into \eqref{c2} we obtain
\begin{equation}\label{c3}
\epsilon (\Re \la Mc,\zeta v\ra +\|\zeta v\|^2 \|c\|^2) \geq
O(\epsilon^2),
\end{equation}
where $O(\epsilon^2)$ is the Landau symbol for denoting a
(polynomial in this case) expression divisible by $\epsilon^2$. If
$c\neq 0$, we can take $\zeta$ of small modulus such that $d:=\Re
\la Mc,\zeta v\ra +\|\zeta v\|^2 \|c\|^2\neq 0$. But then, taking
$|\epsilon|<<1$ such that $\epsilon \cdot d<0$ we contradict
\eqref{c3}. Therefore if $-M$ is not asymptotically stable then
$c=0$. Thus $\v_t(z)=e^{-tM} z$ is a semigroup of linear fractional
maps from $\B^n$ into $\B^n$ (because $-M$ is dissipative) whose
infinitesimal generator is $\tilde{G}$. Hence, also $G$ is an
infinitesimal generator of a semigroup of linear fractional
self-maps of $\B^n$.

Next, we assume $-M$ is asymptotically stable, that is, all its
eigenvalues have negative real part. In particular $M$ and $ M^\ast$
are invertible and there exists a unique $v\in \C^n$ such that
$-M^\ast v=c$. Let $A_t:=\exp(-tM)$ and $c_t:=(\exp(-tM^\ast)-I)v$.
Notice that both $A_t$ and $c_t$ are defined for all $t\in
[0,+\infty)$. Moreover set $\tilde{\v}_t(z)=(A_t z)/(\la z, c_t\ra
+1)$. {\sl A priori} since $\|c_t\|$ can be strictly greater than
$1$ for some $t>0$, $\tilde{\v}_t$ might  not be defined in all
$\B^n$. However, once fixed $t_0\in (0,+\infty)$, there exists
$r=r(t_0)>0$ and $\epsilon=\epsilon(t_0)>0$  such that for all $z\in
\B^n_r:=\{z\in \B^n: \|z\|<r\}$ and $t\in S(t_0,\epsilon):=\{t\in
[0,+\infty):|t-t_0|<\epsilon\}$, the map $\tilde{\v_t}$ is well
defined in $\B^n_r\times S(t_0,\epsilon)$, holomorphic in the first
variable and holomorphic in the second variable. A direct
computation shows that $\tilde{\v}_t(z)$ satisfies
\[
\frac{\de}{\de t}\tilde{\v}_t(z)=\tilde{G}(\tilde{\v}_t(z)), \quad
\tilde{\v}_0(z)=z
\]
in $\B^n_r\times S(t_0,\epsilon)$. Therefore, by uniqueness of
Cauchy-type problem (see, {\sl e.g.}, \cite{Horm}), if $\psi_t$ is
the semigroup of holomorphic self-maps of $\B^n$ whose infinitesimal
generator is $\tilde{G}$ it follows that $\psi_t(z)\equiv
\tilde{\v}_t(z)$ in $\B^n_r\times S(t_0,\epsilon)$. By holomorphic
continuation for any fixed $t\in S(t_0,\epsilon)$ it follows that
$\psi_t(z)=\tilde{\v}_t(z)$ for all $z\in \B^n$. Since $\psi_t$
sends the unit ball into itself, we have that $\Vert c_t\Vert <1$
for all $t>0$, proving that $\psi_t=\tilde{\v}_t$ is actually a
linear fractional self-map of $\B^n$. By the arbitrariety of $t_0$
we find that $\psi_t$ is a linear fractional semigroup of $\B^n$ and
then $\tilde{G}$---and hence $G$---is the infinitesimal generator of
a semigroup of linear fractional self-maps of $\B^n$, as wanted.

Now we are left with the case $G(z)\neq O$ for all $z\in \B^n$. In
this case $\v_t$ has no fixed points in $\B^n$ for all $t>0$ and
there exists a unique common Denjoy-Wolff point that, up to
rotations, we may assume to be $e_1$. We transfer our considerations
to the Siegel half-plane $\H^n:=\{(\zeta,w)\in \C\times \C^n: \Re
\zeta>\|w\|^2\}$ by means of the Cayley transform $C:\B^n\to \H^n$
given by $C(z_1,z''))=(\frac{1+z_1}{1-z_1}, \frac{z''}{1-z_1})$. The
infinitesimal generator $G$ is mapped to
$dC_{C^{-1}(\zeta,w)}(G(C^{-1}(\zeta,w))$. For a vector $v\in \C^n$
we will write $v=(v_1,v'')\in \C\times\C^{n-1}$. A direct
computation shows that
\begin{equation*}
\begin{split}
dC_{C^{-1}(\zeta,w)}(G(C^{-1}(\zeta,w))&=\frac{\zeta+1}{2}\left(
                                                            \begin{array}{cc}
                                                              \zeta+1 & 0 \\
                                                              w & {\sf Id} \\
                                                            \end{array}
                                                          \right)\left(
                                                                   \begin{array}{c}
                                                                     a_1 \\ a''
                                                                   \end{array}
                                                                 \right)\\&-\frac{\la
                                                                 \left(\begin{array}{c}\zeta-1\\2w\end{array}\right),a\ra}{2(\zeta+1)}\left(
                                                                                             \begin{array}{cc}
                                                                                               \zeta+1 & 0 \\
                                                                                               w & {\sf Id} \\
                                                                                             \end{array}
                                                                                           \right)\left(\begin{array}{c}
                                                                     \zeta-1\\
                                                                     2w
                                                                   \end{array}
                                                                 \right)\\
& -\frac{1}{2}\left(
\begin{array}{cc}
                                                              \zeta+1 & 0 \\
                                                              w & {\sf Id} \\
                                                            \end{array}
                                                          \right)A\left(
                                                                    \begin{array}{c}
                                                                      \zeta-1 \\
                                                                      2w \\
                                                                    \end{array}
                                                                  \right)\\&-\frac{\la
                                                                 \left(\begin{array}{c}\zeta-1\\2w\end{array}\right),b\ra}{2(\zeta+1)}\left(
                                                                                             \begin{array}{cc}
                                                                                               \zeta+1 & 0 \\
                                                                                               w & {\sf Id} \\
                                                                                             \end{array}
                                                                                           \right)\left(\begin{array}{c}
                                                                     \zeta-1\\
                                                                     2w
                                                                   \end{array}
                                                                 \right).
\end{split}
\end{equation*}
Therefore $dC_{C^{-1}(\zeta,w)}(G(C^{-1}( \zeta,w)) $ is of the form
$N+Mz+Q\la z,z\ra$ with $N\in \C^n$, $M\in \C^{n\times n}$ and $Q$ a
quadratic form on $\C^n\times\C^n$. We examine the quadratic terms.
As a matter of notation, we write
$A=\left(\begin{smallmatrix} A_{11} & A_{12} \\ A_{21} & A_{22} \\
\end{smallmatrix}\right)$, for $A_{11}\in \C$, $A^t_{12}, A_{21}\in
\C^{n-1}$ and $A_{22}\in \C^{(n-1)\times (n-1)}$. Moreover, if $v,w$
are vector in $\C^m$ for $m\geq 1$ we write $vw:=\sum_{j=1}^m
v_jw_j$. Thus we have
\begin{equation*}
\begin{split}
Q((\zeta,w),(\zeta,w))&=\frac{a_1\zeta}{2}\left(
                                                       \begin{array}{c}
                                                          \zeta \\
                                                         w \\
                                                       \end{array}
                                                     \right)
-\frac{\overline{a_1}\zeta+2w\overline{a''}}{2}\left(
                                                       \begin{array}{c}
                                                          \zeta \\
                                                         w \\
                                                       \end{array}
                                                     \right)
-\frac{1}{2}\left(
              \begin{array}{c}
                A_{11}\zeta^2+2\zeta w A_{12} \\
                A_{11}\zeta w+2(wA_{12})w \\
              \end{array}
            \right)\\
&-\frac{\overline{b_1}\zeta+2w\overline{b''}}{2}\left(
                                                       \begin{array}{c}
                                                          \zeta \\
                                                         w \\
                                                       \end{array}
                                                     \right)
\\&=\left(
\begin{array}{c}
\frac{\zeta^2}{2}(a_1-\overline{a_1}-A_{11}-\overline{b_1})-
\zeta w(\overline{a''}+A_{12}+\overline{b''}) \\
\frac{\zeta w}{2}(a_1-\overline{a_1}-A_{11}-\overline{b_1})
-[w(\overline{a''}+A_{12}+\overline{b''})] w\\
\end{array}
\right).
\end{split}
\end{equation*}
By Lemma \ref{JWC-for-infin}.(1) applied to $G(z)$ we obtain the
following equality:
\[
a_1-\overline{a_1}-A_{11}-\overline{b_1}=0,
\]
while, applying Lemma \ref{JWC-for-infin}.(3) to $G(z)$ we obtain
\[
\overline{a''}+A_{12}+\overline{b''}=O,
\]
and therefore $Q\equiv 0$.

Thus $\tilde{G}(z)=Mz+N$ for some $n\times n$ matrix and some vector
$N\in \C^n$. Let $\psi_t$ be the semigroup of holomorphic self-maps
of $\H^n$ of whom $\tilde{G}$ is the infinitesimal generator. Let
\[
\tilde{\v}_t(z):=e^{tM}z+\left(\int_0^t e^{sM}ds \right)N.
\]
Then $\tilde{\v}_t$ is defined for all $z\in \C^n$ and all $t\in
[0,+\infty)$. A direct computation shows that
\[
\frac{\de}{\de t} \tilde{\v}_t(z)=\tilde{G}(\tilde{\v}_t(z))
\]
for all $z$ and $t$. Thus by the uniqueness of solution of partial
differential equations, it follows that $\tilde{\v}_t=\psi_t$ for
all $t$ and in particular the semigroup associated to $\tilde{G}$ is
linear, which, going back to the ball, means that $G$ is the
infinitesimal generator of a linear fractional semigroup of $\B^n$
as wanted.
\end{proof}

\begin{corollary} Let $(\v_t)$ be a semigroup of holomorphic self-maps of
$\B^n$ with infinitesimal generator   $G\in Hol(\B^n,\C^n)$. The
following are equivalent.
\begin{enumerate}
\item $(\v_t)$ is a group of holomorphic self-maps of
$\B^n$.
\item There exist $a\in \C^n$ and $A\in \C^{n\times n}$ (not
all of them zero) such that $G(z)=a-\la z,a\ra z-Az$ with
\begin{equation}
\Re \la Az, z\ra = 0,
\end{equation}
for all $z\in \C^n$.
\end{enumerate}
\end{corollary}
\begin{proof}
If $(\v_t)$ is a group of holomorphic self-maps of $\B^n$, then the
family of functions $(\psi _t) :=(\v_{-t})$ is a semigroup of linear
fractional maps of $\B^n$ and its infinitesimal generator is $-G$.
Therefore, bearing in mind that $G$ and $-G$ are infinitesimal
generator of linear fractional maps,
 the above theorem implies that we have
\begin{equation}
\Re \la Au, u\ra  = 0
\end{equation}
for all $u$ and \eqref{constrain} implies that $b=0$. Conversely, if
$G(z)=a-\la z,a\ra z-Az$ with $\Re \la Az, z\ra=0$ for all $z\in
\C^n$,  the above theorem shows that both $G$ and $-G$ are
infinitesimal generators and the Cauchy problems
\begin{equation*}
\frac{\partial g}{\partial t}=\pm G\circ g,\  \hbox{ with} \ g(0)=z
\end{equation*}
have solutions in $[0,+\infty)$ for all $z\in\B^n$. Therefore $G$ is
the infinitesimal generator of a group of automorphisms.
\end{proof}

\section{Infinitesimal generators in dimension one and Koenigs Functions}

Koenigs functions and infinitesimal generators are related in a very
concrete way. An explicit reference for the following result is not
really available but essentially the idea of the proof is given in
\cite{Siskakis-tesis}.

\begin{proposition}
\label{Unival-VectorField} Let $(\varphi _{t})$ be a non-trivial
semigroup in $\mathbb{D}$ with infinitesimal generator $G$   and
Koenigs function $h$.
\begin{enumerate}
\item If $(\varphi _{t})$ has Denjoy-Wolff point $\tau \in \mathbb{D}$ then  $h$ is the unique
holomorphic function from $\mathbb{D}$ into $\mathbb{C}$ such that
\begin{enumerate}
\item[(i)] $h^{\prime }(z)\neq 0,$ for every $z\in \mathbb{D},$
\item[(ii)] $h(\tau )=0$ and $h^{\prime }(\tau )=1,$
\item[(iii)] $h^{\prime }(z)G(z)=G^{\prime }(\tau )h(z),$ for
every $z\in \mathbb{D}.$
\end{enumerate}
\item If $(\varphi _{t})$ has Denjoy-Wolff point $\tau
\in \partial \mathbb{D}$ then $h$ is the unique holomorphic function
from $\mathbb{D}$ into $\mathbb{C}$ such  that:
\begin{enumerate}
\item[(i)] $h(0)=0,$
\item[(ii)] $h^{\prime }(z)G(z)=1,$ for every $z\in \mathbb{D}.$
\end{enumerate}
\end{enumerate}
\end{proposition}

For further reference  we state here the following  simple fact.

\begin{lemma}
\label{lemma mz+n}Let $p(z)=mz+n$ be a complex polynomial. Then  $\func{Re}%
p(z)\geq 0$   for every $z\in \mathbb{D}$ if and only if
$\func{Re}(n)\geq |m|.$
\end{lemma}

Our first result is a translation to one variable of the result
contained in the previous section. This result extends and, somehow,
completes \cite[Proposition 3.5.1]{Shoiket}.

\begin{theorem}
\label{Tma Basico}Let $G:\mathbb{D\rightarrow C}$ be a  holomorphic
function. Then, the following are equivalent:

\begin{enumerate}
\item The map $G$ is the infinitesimal generator of a semigroup of $\mathrm{LFM}(\mathbb{D},\mathbb{D}).$

\item The map $G$ is a polynomial of degree at most two and
satisfies the following boundary flow condition
\begin{equation*}
\func{Re}(G(z)\overline{z})\leq 0,\text{ for all }z\in \partial
\mathbb{D}.
\end{equation*}

\item The map $G$ is a polynomial of the form $G(z)=\alpha z^{2}+\beta z+\gamma$ with $\func{%
Re}(\beta)+|\alpha+\overline{\gamma}|\leq 0.$

\item The map $G$ is a polynomial of the form
\begin{equation*}
G(z)=a-\overline{a}z^{2}-z(mz+n),\text{ }z\in \mathbb{D}
\end{equation*}%
with $a,m,n\in \mathbb{C}$ and $\func{Re}(n)\geq |m|.$

\item The map $G$ is the infinitesimal generator  of a semigroup in $\mathbb{D}$
and it is a polynomial of degree at most two.
\end{enumerate}

Moreover, $\func{Re}(n)=m=0$ in statement $(4)$ if and only if $\func{Re}%
(\beta)=|\alpha+\overline{\gamma}|=0$ in statement $(3)$ if and
only if equality holds for all $z\in \partial \mathbb{D}$ in
statement $(2)$ if and only if $G$ is the infinitesimal generator
of a semigroup of $\Aut(\mathbb{D})$.
\end{theorem}

\begin{proof}
By Theorem \ref{main}, statement (1) is equivalent to (4).
Statement (4) is then equivalent to (2), (3) and (5) by direct
computations using Lemma \ref{lemma mz+n}. Finally, the last
assertion follows from the fact that if $G$ is the infinitesimal
generator of a semigroup in $\D$, then this semigroup is composed
of automorphisms of $\D$ if and only if $-G$ is the infinitesimal
generator of a semigroup in $\D$ as well.
\end{proof}

This theorem clearly implies that not every polynomial of degree two
can be realized as an infinitesimal generator of a semigroup in
$\mathbb{D}.$ At the same time, it also suggests to analyze
carefully those complex polynomials $G$ of degree zero, one, and two
which are infinitesimal generators of  semigroups in $\mathbb{D}$.

Trivially, if $G$ is of degree zero, then $G$ is an infinitesimal
generator of a semigroup in $\mathbb{D}$ if and only if the constant
is zero and, in this case, the corresponding semigroup is the
trivial one. The next two propositions describe what happens when
the degree is one and two.

\begin{proposition}\label{poldeg1}
Let $G(z)=\lambda (z-c)$, $\lambda\neq 0$, be a complex polynomial
of degree one. Then
 $G$ is the infinitesimal generator of a semigroup in
$\mathbb{D}$ (which is necessarily a semigroup of
$\mathrm{LFM}(\mathbb{D},\mathbb{D})$) if and only if
$\func{Re}\lambda +\left\vert \lambda c\right\vert \leq 0.$
Moreover, if $G$ is an infinitesimal generator the associated
semigroup is  given by
\begin{equation*}
\varphi _{t}(z)=e^{\lambda t}z+c(1-e^{\lambda t}),\qquad z\in
\mathbb{D},
\end{equation*}
and $c\in\overline{\D}$   is its Denjoy-Wolff point. Furthermore,

\begin{enumerate}
\item if $\func{Re}\lambda  =0$ then $c=0$ and the semigroup is a
neutral-elliptic semigroup  of $\Aut(\mathbb{D})$.

\item if $0<|\Re\lambda|<|\lambda|$ then $|c|<1$, $\left\vert \func{Re}
\lambda \right\vert \geq |\func{Im}\lambda
|\dfrac{|c|}{\sqrt{1-|c|^{2}}}$ with $\func{Re}\lambda <0$ and the
semigroup is attractive-elliptic.

\item if $|\Re \lambda|=|\lambda|$ then $|c|=1$, $\lambda \in (-\infty
,0)$ and the semigroup is hyperbolic.
\end{enumerate}
\end{proposition}

\begin{proof}
Theorem \ref{Tma Basico} implies that $G(z)=\lambda (z-c)$ is an
infinitesimal generator of a semigroup in $\mathbb{D}$ (which is
necessarily a semigroup of linear fractional maps) if and only if
$\func{Re}\lambda +\left\vert \lambda c\right\vert \leq 0.$ Note
that the latter inequality   implies that $c\in
\overline{\mathbb{D}}$. Moreover, $G(c)=0$ and thus in case $G$ is
the infinitesimal generator of a semigroup in $\D$, the point $c$
is the Denjoy-Wolff point of the associated semigroup (this
follows for instance from Berkson-Porta's Theorem~\ref{24} which
shows that, unless $G\equiv 0$ then it has a unique zero in $\D$
which is exactly the Wolff-Denjoy point of the associated
semigroup).

Now assume that $G$ is the infinitesimal generator of the
semigroup $(\v_t)$ in $\D$. If $|\Re\lambda|<|\lambda|$, then the
inequality $\func{Re}\lambda +\left\vert \lambda c\right\vert \leq
0$ implies $c\in \mathbb{D}$. By Proposition 2.1, the Koenigs
function $h$ of the semigroup satisfies
\begin{equation*}
h^{\prime }(z)(z-c)=h(z),\qquad z\in \mathbb{D},
\end{equation*}
$h(c)=0$ and $h^{\prime }(c)=1.$ That is, $h(z)=z-c$ and $\varphi
_{t}(z)=h^{-1}(e^{G^{\prime }(c)t}h(z))=e^{\lambda
t}z+c(1-e^{\lambda t})$ for all $z.$ Now (1) and (2) follow easily
from direct computations.

If $|\Re\lambda|=|\lambda|$ then $\lambda=\Re\lambda<0$ and $c\in
\partial \mathbb{D}$ again by the inequality $\func{Re}\lambda +\left\vert \lambda c\right\vert \leq
0$ and  the semigroup is hyperbolic.
\end{proof}

The previous proposition implicitly says that linear infinitesimal
generators of semigroups in $\mathbb{D}$ are in one-to-one
correspondence with affine semigroups. In other words, for seeing
non-linear phenomena, we have to deal with polynomials of degree two
or more.

\begin{proposition}\label{poldeg2}
Let $G(z)=\lambda (z-c_{1})(z-c_{2})$ be a complex polynomial of
degree two with $\lambda =|\lambda |e^{i\theta }$, $\lambda\neq 0$.
Then, $G$ is the infinitesimal generator of a semigroup in
$\mathbb{D}$ (necessarily a semigroup of
$\mathrm{LFM}(\mathbb{D},\mathbb{D})$) if and only if
\begin{equation}\label{p1r}
\func{Re}(e^{i\theta }c_{1}+e^{i\theta }c_{2})\geq 0\text{ and }%
(|c_{1}|^{2}-1)(1-|c_{2}|^{2})\geq [\func{Im}(e^{i\theta
}c_{1}-e^{i\theta }c_{2})]^{2}.
\end{equation}
Moreover if $G$ is the infinitesimal generator of a semigroup in
$\mathbb{D}$ then $c_1\in\overline{\D}$ is the Denjoy-Wolff point of
$(\v_t)$ and the following   are the only possible cases:

\begin{enumerate}
\item if $c_1=c_2\in \partial \mathbb{D}$ then
$\func{Re}(e^{i\theta }c_1)\geq 0$ and $(\v_t)$ is a parabolic
semigroup. Moreover in this case, $(\v_t)$ is a parabolic
semigroup of $\Aut(\mathbb{D})$ if and only if
$\func{Re}(e^{i\theta }c_1)=0.$

\item if $c_1 \in \de\D$ and $c_2\in \C\setminus (\D\cup\{c_1\})$
then $e^{i\theta }(c_{2}-c_{1})\in (0,+\infty )$ and $(\v_t)$ is a
hyperbolic semigroup. Moreover in this case, $(\v_t)$ is a
hyperbolic semigroup of $\Aut(\mathbb{D})$ if and only if $c_2\in
\de \D$ if and only if $e^{i2\theta }c_{1}c_{2}=-1$.

\item if $c_1\in\D$ then $c_2\in\C\setminus\D$ and $(\v_t)$ is an
elliptic semigroup. Moreover:
\begin{itemize}
\item[a)] if $c_{2}\in \partial \mathbb{D}$ then the semigroup is
attractive-elliptic with two fixed points in
$\overline{\mathbb{D}}$ and  $e^{i\theta }(c_{2}-c_{1})\in
(0,+\infty )$. \item[b)] if $c_2\in \C\setminus\overline{\D}$ and
$c_2\overline{c_{1}}\neq 1$ then the semigroup is
attractive-elliptic   with only one fixed point in
$\overline{\mathbb{D}}$ and $\func{Re}(e^{i\theta
}(c_{1}+c_{2}))\in (0,+\infty )$, $\func{Im}(e^{i\theta
}(c_{1}-c_{2}))\in \lbrack -\beta ,\beta ],$ where $\beta
:=\sqrt{(|c_{1}|^{2}-1)(1-|c_{2}|^{2})}>0.$ \item [c)] if $c_1\neq
0$ and $c_2\overline{c_{1}}=1$ then the semigroup is a
neutral-elliptic semigroup of $\Aut(\mathbb{D})$ and
$\func{Re}(e^{i\theta}c_1)=0$.
\end{itemize}
\end{enumerate}
\end{proposition}

\begin{proof}
First of all, a direct computation from  Theorem \ref{Tma Basico}
 shows that $G$ is the infinitesimal generator of a semigroup of
linear fractional maps of $\D$ if and only if \eqref{p1r} holds.

Now assume that $G$ is the infinitesimal generator of a semigroup
$(\v_t)$ of $\mathrm{LFM}(\mathbb{D},\mathbb{D})$ and then
\eqref{p1r} holds. We first recall that (as a consequence of
Berkson-Porta's Theorem~\ref{24}) either $c_1$ or $c_2$ is the
Denjoy-Wolff point of $(\v_t)$. Without loss of generality we can
assume that $c_1$ is the Denjoy-Wolff point of $(\v_t)$. From
Theorem~\ref{24} one sees that $c_2\not\in\D$. Furthermore, as a
consequence of the Schwarz lemma or the Julia-Wolff-Caratheodory
theorem (see \cite[Theorem 1]{Contreras-Diaz-Pom-SCAND} or
\cite{Shoiket}) one can show that $\Re G'(c_1)\leq 0$ and, if
$c_1\in\de\D$ then $\Re G'(c_1)=0$ if and only if $(\v_t)$ is a
parabolic semigroup of linear fractional maps. Taking these into
account, direct computations from \eqref{p1r} give statements (1),
(2), and (3), with the possible exception of the characterization
of semigroups of automorphisms.

In order to obtain the characterization of semigroup of
automorphisms, it is enough to realize that $(\v_t)$ is a semigroup
of automorphisms if and only if $-G$ is an infinitesimal generator
for a semigroup in $\D$ and apply \eqref{p1r} to $-G$.
\end{proof}

\begin{remark}
Assume $G$ is as in Proposition \ref{poldeg2}.(3).(c) and conjugate
the semigroup $(\v_t)$ with the automorphism
$T(z)=(c_1-z)(1-\overline{c_1}z)^{-1}$. The new semigroup $(T\circ
\v_t\circ T)$ has infinitesimal generator given by
$\tilde{G}(z)=T'(T(z))\cdot G(T(z))$. A direct computation shows
that $\tilde{G}(z)=\lambda (\overline{c_1})^{-1}(1-|c_1|^2)z$ and
thus we are in the case of Proposition~\ref{poldeg1}.(1).
\end{remark}

Since an infinitesimal generator $G$ generates a semigroup of
automorphisms of $\D$ if and only if $-G$ is an infinitesimal
generator of a semigroup in $\D$, in the Berkson-Porta
representation of an infinitesimal generator  $G(z)=(z-\tau
)(\overline{\tau }z-1)p(z)$, parabolic semigroups of
$\Aut(\mathbb{D})$ appear exactly when $\tau \in
\partial \mathbb{D}$ and $p=i\beta $ for some real $\beta \neq 0$ and
elliptic (necessary neutral-elliptic) semigroups of $\Aut(\mathbb{D}%
) $ are exactly generated when $\tau \in \mathbb{D}$ and $p=i\beta $
for some real $\beta \neq 0$. Using Theorem 2.3, we  can also
explain the hyperbolic case.

Before giving the statement, we recall that if $\H:=\{w\in \C:\Re
w>0\}$ is the right half plane, a {\sl Cayley transform} with pole
$\tau\in\de\D$ is any biholomorphic map $C:\D\to\H$ such that
$\lim_{\D\ni z\to \tau}|C(z)|=\infty$. It is well known that every
Cayley transform is a linear fractional map.

\begin{proposition}
Let $G:\D\to\C$ be a holomorphic function. Then $G$ is the
infinitesimal generator of a hyperbolic semigroup of
$\Aut(\mathbb{D})$ if and only if $\,G(z)=(z-\tau )(\overline{\tau
}z-1)p(z),$ with $\tau \in \partial \mathbb{D}$ and $p$ is a
Cayley transform with   pole  $\tau$.
\end{proposition}

\begin{proof} $(\Rightarrow )$ Since $G$ generates
a semigroup $(\v_t)$ of hyperbolic automorphisms then there exists
$\sigma\in\de\D\setminus\{\tau\}$ such that $G(\sigma)=0$ (such a
point $\sigma$ is the second common fixed point for $(\v_t)$).
According to Theorem \ref{Tma Basico} we have then
$G(z)=\lambda(z-\tau)(z-\sigma)$. From this we get
$p(z)=\lambda(z-\sigma)(\overline{\tau}z-1)^{-1}$. Now $p(\sigma)=0$
and since linear fractional maps on the Riemann sphere maps circles
onto circles, we see that $p(\de \D)=\de \H$ proving that $p$ is a
Cayley transform with pole $\tau$.

$(\Leftarrow )$ Using Theorem \ref{Tma Basico}, we   see that $G$ is
the infinitesimal generator of a non-trivial semigroup of
$\mathrm{LFM}(\mathbb{D},\mathbb{D}).$ It is enough to show that the
semigroup is composed of automorphisms of $\D$. Indeed $p$ is not
constant, and then $G$ cannot generate parabolic or elliptic groups.
To this aim, using again Theorem \ref{Tma Basico} and writing $p(z)=$ $\dfrac{az+b}{%
\overline{\tau }z-1}$, we have only  to check that $\overline{a}=\tau b$ and $%
\func{Re}(b-a\tau )=0.$ Since $p$ is bijective, $p(\partial
\mathbb{D})=\partial \mathbb{H}$ and, in particular,
$\func{Re}p(-\tau )=0$ thus $\func{Re}( a(-\tau )+b)=0.$ Therefore,
it only remains to prove the other condition. For every $z\in
\partial \mathbb{D},$ we have
\begin{equation*}
0=\func{Re}((az+b)(1-\overline{z}\tau ))=\func{Re}(b-a\tau )+\func{Re}((%
\overline{a}-b\tau )\overline{z})=\func{Re}((\overline{a}-b\tau )\overline{z}%
).
\end{equation*}%
Therefore $|\overline{a}-b\tau |=0$ as needed.
\end{proof}

We end this section looking at   relationships between Koenigs
functions and semigroups of linear fractional maps. As customary,
$\infty ^{-1}$ means $0$.

\begin{proposition}\label{kperlfm}
Let $h:\mathbb{D\rightarrow C}$ be a  holomorphic function.

\begin{enumerate}
\item The map $h$ is the Koenigs function associated to a non
trivial-elliptic semigroup of
$\mathrm{LFM}(\mathbb{D},\mathbb{D})$ if and only if
\begin{equation*}
h(z)=(1-\beta ^{-1}\tau )\frac{z-\tau }{1-\beta ^{-1}z}
\end{equation*}%
for some $\tau \in \mathbb{D}$ and $\beta \in \mathbb{C}_{\infty
}\setminus \mathbb{D}$.
\item The map $h$ is the Koenigs function associated to a
hyperbolic semigroup of $\mathrm{LFM}(\mathbb{D},\mathbb{D})$ if
and only if there exists $\alpha \in (-\infty ,0)$ such that
\begin{equation*}
h(z)=\alpha \log \left( \frac{1-\overline{\tau }z}{1-\beta
^{-1}z}\right)
\end{equation*}%
where $\tau \in \partial \mathbb{D},$ $\beta \in
\mathbb{C}_{\infty }\setminus \mathbb{D\,\ }$and $\log $ denotes
the principal branch of the logarithm.

\item The map $h$ is the Koenigs function associated to a
parabolic semigroup of $\mathrm{LFM}(\mathbb{D},\mathbb{D})$ if
and only if there exists $\alpha \neq 0$ with $\func{Re}\alpha
\leq 0$ such that
\begin{equation*}
h(z)=\alpha \frac{z}{z-\tau }
\end{equation*}%
where $\tau \in \partial \mathbb{D}$.
\end{enumerate}
\end{proposition}

\begin{proof}
(1) Let $\eta_\tau(z):=(\tau-z)(1-\overline{\tau}z)^{-1}$. Notice
that $\eta_\tau\in\Aut(\D)$ and $\eta_\tau^{-1}=\eta_\tau$.
Conjugating $(\v_t)$ with $\eta_\tau$ we obtain a semigroup
$(\phi_t)$ with Wolff-Denjoy point $0$. By uniqueness of the Koenigs
function, if $h$ is the Koenigs function of $(\v_t)$ and $g$ is the
Koenigs function of $(\phi_t)$, we have
\begin{equation}\label{we-de}
h(z)=(|\tau|^2-1)g(\eta_\tau(z)).
\end{equation}
Since the infinitesimal generator $F$ of $(\phi_t)$ is given by
$F(z)=\eta_\tau'(z) G(\eta_t(z))$ and, by Theorem \ref{Tma
Basico}, it is of the form $-Mz^2-Nz$ (note that $F(0)=0$ because
$\phi_t(0)=0$ for all $t$). Now, a direct computation using
Proposition \ref{Unival-VectorField}.(1.iii) shows that $g$ is of
the claimed form and, by \eqref{we-de}, so is $h$.

(2) Let $G$ be the infinitesimal generator of $(\v_t)$; and let
$\tau\in\de \D$ be the Denjoy-Wolff point of $(\v_t)$. Using
 Propositions \ref{Unival-VectorField}.(2.iii), \ref{poldeg1}.(3), and
\ref{poldeg2}.(2) a direct computation yields that the Koenigs
function $h$ is of the form
\begin{equation*}
\frac{1}{G^{\prime }(\tau )}\left( \log (1-\overline{\tau }z)-\log
(1-\frac{1}{\beta }z)\right) ,
\end{equation*}%
where  $\beta \in \mathbb{C}_{\infty }\setminus \mathbb{D}$ and
$\log $ is the principal branch of the logarithm. Now, we recall
that $G^{\prime }(\tau )\in(-\infty ,0)$
\cite{Contreras-Diaz-Pom-SCAND} and since $1-\overline{\tau
}z,1-\beta ^{-1}z\in \mathbb{H},$ for all $z\in \mathbb{D}$,
\begin{equation*}
\limfunc{Arg}(1-\overline{\tau }z)-\limfunc{Arg}(1-\frac{1}{\beta }z)=%
\limfunc{Arg}\left( \frac{1-\overline{\tau }z}{1-\beta
^{-1}z}\right) ,
\end{equation*}%
where $\limfunc{Arg}$ denotes the principal argument.

(3) Once more, a direct computation from Proposition
\ref{poldeg2}.(1) and Proposition \ref{Unival-VectorField}.(2.iii),
implies that Koenigs functions associated to parabolic semigroups of
$\mathrm{LFM}(\mathbb{D},\mathbb{D})$ are exactly those of the form
$\alpha \dfrac{z}{z-\tau },$ with $\tau \in
\partial \mathbb{D}$  the Denjoy-Wolff of the semigroup and
$\alpha $   a given number which, by Proposition
\ref{Unival-VectorField}.(2), is $\alpha \neq 0$ and
$\func{Re}\alpha \leq 0.$

Finally, we point out that proving the converse in each of the
three cases is just a direct (and lengthy) computation.
\end{proof}

It is clear from the proof that the point $\tau $ in the statement
of the previous theorem is the Denjoy-Wolff point of the
corresponding semigroup, while the point $\beta $ in cases $1$ and
$2$ is exactly the repulsive fixed point
 of each iterate, seen as a M\"obius transform of the Riemann sphere.

\section{The Embedding Problem}

The abstract embedding problem is a rather classical one and has
been treated from the times of Abel. In its most general form the
embedding problem can be stated as follows:  given a space $X$ in
some category (topological, differential, holomorphic) and a map
$f:X\rightarrow X$ in the same category,  determine whether it is
possible to construct a   semigroup $(f_{t})$
over $X$, continuous in $t$ with iterates $f_t$ in the same category of $f$,
 such that $f=f_{1}$.

When $X=\mathbb{D}$ and $f$ is a linear fractional map, the
embedding problem has been treated by several authors (see
\cite{Khat-Shoikhet} and references therein) and   it is known that,
in general, it has a negative answer.

As stated in the introduction, in this paper we prove that a linear
fractional self-map of $\D$ can be embedded in a semigroup of $\D$
if and only if it can be embedded into a semigroup of linear
fractional self-maps of $\D$, obtaining a precise characterization
of those linear fractional map which can be embedded.

 The proof of
this result requires several tools, most of them based on model
theory and some of them quite recent. In fact, we will need a new
representation, with strong uniqueness, for hyperbolic semigroups in
$\mathbb{D}$ which could have some interest in its own. Such a model
is a continuous version of the classical Valiron's construction,
(see, {\sl e.g.} \cite{Bracci-Poggi}).

\begin{proposition}\label{intertwining}
Let $(\varphi _{t})$ be a hyperbolic semigroup in $\mathbb{D}$
with Denjoy-Wolff point $\tau \in \partial \mathbb{D}$ and
associated Koenigs
function $h$. Then there exists a univalent function $\sigma :\mathbb{%
D\rightarrow H}$ $\ $such that $|\sigma (0)|=1$ and%
\begin{equation*}
\sigma \circ \varphi _{t}=\varphi _{t}^{\prime }(\tau )\sigma
,\text{ for every }t\geq 0.
\end{equation*}%
Indeed, there exists $\alpha \in (-\infty ,0)$ such that $\sigma
=\sigma (0)e^{\alpha h}.$

Moreover, if a non-constant holomorphic (a priori, not necessary
univalent) function $\rho \in \mathrm{Hol}(\mathbb{D},\mathbb{H})$
satisfies $|\rho (0)|=1$ and $\rho \circ \varphi _{t}=\varphi
_{t}^{\prime }(\tau )\rho ,$ for every $t\geq 0,$ then $\rho
=\sigma .$
\end{proposition}

\begin{proof}
Since every iterate $\varphi _{t}$ $(t>0)$ is hyperbolic, we can
consider the Valiron  normalization (see, {\sl e.g.}
\cite{Bracci-Poggi}) with respect to   $0$ and obtain a non-constant
holomorphic function $\sigma _{t}\in
\mathrm{Hol}(\mathbb{D},\mathbb{H})$ such that $|\sigma _{t}(0)|=1$
and
\begin{equation*}
\sigma _{t}\circ \varphi _{t}=\varphi _{t}^{\prime }(\tau )\sigma
_{t}.
\end{equation*}%
Such a map $\sigma_t$  is univalent for all $t\geq 0$ because $\v_t$
is. Now recall that by the chain-rule for non-tangential
derivatives, $ \lbrack 0,+\infty )\ni t\mapsto \varphi _{t}^{\prime
}(\tau )\in \lbrack 0,1]$ is a measurable algebraic homomorphism
from $(\R,+)$ and $(\R^\ast,\cdot)$. Therefore, there exists $\alpha
\in (-\infty ,0)$ such that $\varphi _{t}^{\prime }(\tau )=e^{\alpha
t}.$

Set $\sigma =\sigma _{1}.$ Then, for every $n\in \mathbb{N},$ we
have $\sigma \circ \varphi _{n}=e^{n\alpha t}\sigma =\varphi
_{n}^{\prime }(\tau )\sigma$. Hence, both    $\sigma ,$ $\sigma
_{n}$ are intertwining functions for $\v_n$ and, according to the
strong uniqueness result proved in \cite[Proposition
6]{Bracci-Poggi}, we deduce that there exists $c_{n}>0$ such that
$\sigma =c_{n}\sigma _{n}.$ From the choice of our normalization,
$|\sigma (0)|=|\sigma _{n}(0)|=1$, and so $\sigma =\sigma _{n}.$ A
similar argument also shows that $\sigma =\sigma _{t}$ for every
positive and rational $t$ and, finally,   continuity in $t$ of the
semigroup implies that $\sigma =\sigma _{t}$ for every $t>0$.

Since $\sigma (\mathbb{D})\subset \mathbb{H},$ we can consider $\log
\sigma (z),$ where $\log $ denotes the principal branch of the
logarithm and, by definition, $\log (\sigma  \circ \varphi
_{t})=\log \sigma +\alpha t$ for all $t\geq0$. Therefore, the
function $\widehat{\sigma }:=\log \sigma (z)-\log \sigma
(0)\in \mathrm{Hol}(\mathbb{D},\mathbb{C})$ satisfies   $\widehat{\sigma }%
(0)=0$ and $\widehat{\sigma }$ $\circ \varphi _{t}=\widehat{\sigma
}+\alpha t.$ Differentiating with respect to $t$ and evaluating at
$t=0,$ we find that
\begin{equation*}
\widehat{\sigma }^{\prime }(z)G(z)=\alpha ,\text{ for every
\thinspace }z\in \mathbb{D},
\end{equation*}%
where $G$ denotes the infinitesimal generator of $(\varphi _{t})$.
According to Proposition \ref{Unival-VectorField}, we conclude that $\dfrac{1%
}{\alpha }\widehat{\sigma }=h.$ Thus $\sigma =\sigma (0)e^{\alpha
h}.$

Finally, the uniqueness assertion follows from the corresponding one
for intertwining mappings. Indeed, if $\rho \in
\mathrm{Hol}(\mathbb{D},\mathbb{H})$ is non constant and satisfies
$\rho \circ \varphi _{t}=\varphi _{t}^{\prime }(\tau )\rho$ for all
$t\geq 0$, then in particular $\rho \circ \varphi _{1}=\varphi
_{1}^{\prime }(\tau )\rho$ and according to \cite[Proposition
6]{Bracci-Poggi} then $\rho=c\sigma$ for some $c>0$. The further
normalization $|\rho (0)|=1$ implies that $\rho=\sigma$.
\end{proof}

Now we can prove the following rigidity result:

\begin{theorem}\label{rigido}
Let $(\v_t)$ be a semigroup in $\D$. If for some $t_0>0$ the iterate
$\v_{t_0}$ is a linear fractional self-map of $\D$ then $\v_t$ is a
linear fractional self-map of $\D$ for all $t\geq 0$.
\end{theorem}
\begin{proof} Up to rescaling we can assume that $t_0=1$. First of all, if $\v_{1}\in\Aut(\D)$ then
$\v_t\in\Aut(\D)$ for all $t\geq 0$ and the result holds. Thus we
assume $\v_{1}$ is not surjective and we consider the possible
dynamical types of $\v_{1}$:

(Attractive-Elliptic case) Since composition of linear fractional
maps is linear fractional, up to conjugation with a suitable
automorphism of $\D$, we can assume that $(\v_t)$ has Denjoy-Wolff
point $0$. Thus
\[
\v_1(z)=\frac{az}{cz+1}, \quad z\in\D
\]
for $a=\v_1'(0)$  (thus $0<|a|<1$) and $c\in\C$. Hence for all
$n\in\mathbb N$
\begin{equation}\label{iteraellip}
\v_n(z)=\frac{a^nz}{c\frac{1-a^n}{1-a}z+1}.
\end{equation}
Now let   $\sigma:\D\to \C$ be the Schr\"oder function of $\v_1$
(see, {\sl e.g.}, \cite{Shapiro-libro}). This is the unique
univalent (because $\v$ is) holomorphic function such that
$\sigma(0)=0$, $\sigma'(0)=1$ and $\sigma\circ \v_1=\v_1'(0)\sigma$;
it is defined as the limit of the sequence $\{\v_n/\v_n'(0)\}$. Thus
a direct computation from \eqref{iteraellip} shows that $\sigma$ is
the linear fractional map given by
\begin{equation}\label{hK}
z\mapsto \frac{z}{\frac{c}{1-a}z+1}.
\end{equation}
 If $h:\D\to\C$ is the Koenigs function of the semigroup $(\v_t)$
then $h(0)=1$, $h'(0)=1$ and $h\circ\v_t=\v_t'(0)h$ for all $t\geq
0$. In particular $h\circ\v_1=\v_1'(0)h$ and by uniqueness
$h=\sigma$, which implies $\v_t(z)=h^{-1}(\v_t'(0)h(z))\in
\mathrm{LFM}(\mathbb{D},\mathbb{D})$ for all $t\geq 0$.

(Hyperbolic case) The semigroup $(\v_t)$ is hyperbolic with
Denjoy-Wolff point $\tau\in\de \D$. Let $\sigma:\D\to\H$ be the
function defined in Proposition \ref{intertwining}. By the very
construction and uniqueness of the Valiron intertwining function
(see \cite{Bracci-Poggi}) it follows that
\begin{equation}\label{valiron}
\sigma (z)=\lim_{\mathbb N\ni n\to \infty}\frac{1-\overline{\tau }\varphi _{n}(z)}{\left\vert 1-%
\overline{\tau }\varphi _{n}(0)\right\vert },\text{ }z\in
\mathbb{D}.
\end{equation}
Now, $\varphi _{n}\in \mathrm{LFM}(\mathbb{D},\mathbb{D})$ for all
$n\in \mathbb{N}$, so the Schwarzian derivates $S_{\varphi
_{n}}\equiv 0$ in $\mathbb{D}$  for all $n\in \mathbb{N}$. Since the
limit in \eqref{valiron} holds uniformly on compacta in
$\mathbb{D}$, then   $S_{\sigma }\equiv 0$ in $\mathbb{D}$ implying
that   $ \sigma \in \mathrm{LFM}(\mathbb{D},\mathbb{H})$. Finally
$\varphi _{t}(z)=\sigma ^{-1}(\v_t'(\tau)\sigma (z))\in
\mathrm{LFM}(\mathbb{D},\mathbb{D})$  for every $t\geq 0$.

(Parabolic case) The semigroup $(\v_t)$ is   parabolic with
Denjoy-Wolff point $\tau \in \partial \mathbb{D}$ and $\varphi
_{t}^{\prime }(\tau )=1$. Let $h$ be the associated Koenigs
function. Thus, $h\in \mathrm{Hol}(\mathbb{D},\mathbb{C})$ is
univalent, verifies $h(0)=0$ and $h\circ \varphi _{t}=h+t$. Since
$\varphi_1 \in \mathrm{LFM}(\mathbb{D},\mathbb{D})\setminus
\Aut(\mathbb{D})$ and it is parabolic,
\begin{equation*}
\lim_{\mathbb N\ni n\to\infty}k_\D(\varphi _{n}(z),\varphi
_{n+1}(z))=0,
\end{equation*}
for every $z\in \mathbb{D}$, where $k_\D$ denotes the hyperbolic
metric in $\mathbb{D}$.  Thus we can apply Baker-Pommerenke's
normalization (see \cite{BakerPomm79}) and obtain a univalent map
$\sigma \in \mathrm{Hol}(\mathbb{D},\mathbb{C})$, which is a uniform
limit on compacta of $\D$ of linear fractional combinations of
iterates of $\v_1$ such that $\sigma (0)=1$ and $\sigma \circ
\varphi =\sigma +1$. Arguing as before,
$\sigma\in\mathrm{LFM}(\mathbb{D},\mathbb{C})$. By \cite[Theorem
3.1]{Con-DiazPom-Abel} there exists $\lambda \in \mathbb{C}$ such
that $h=\sigma +\lambda \in \mathrm{LFM}(\mathbb{D},\mathbb{C}).$
Thus $\varphi _{t}(z)=h^{-1}(h(z)+t)\in
\mathrm{LFM}(\mathbb{D},\mathbb{D})$, for every $t\geq 0$.
\end{proof}

Now the following result is a direct consequence of Theorem
\ref{rigido} and the fact that the automorphisms of $\C_\infty$ are
exactly linear fractional maps.

\begin{theorem}\label{embedding}
Let $\varphi\in\mathrm{LFM}(\mathbb{D},\mathbb{D}).$ Then, the
following are equivalent:

\begin{enumerate}
\item The map $\varphi $ can be embedded into a semigroup in
$\mathbb{D}.$

\item The map $\varphi $ can be embedded into a semigroup of $\mathrm{LFM}(%
\mathbb{D},\mathbb{D}).$

\item The map $\varphi$, thought of  as an element of $\Aut(\C_\infty)$, can be embedded  into a group
$(\varphi _{t})_{t\in \mathbb{R}}$ of $\Aut(\mathbb{C}_{\infty})$
with the property that $\varphi _{t}(\D)\subseteq \D$ for all $t\geq
0$.
\end{enumerate}
\end{theorem}

In the last result of this section, we provide a purely dynamical
analysis of the embedding problem for linear fractional self-maps
of the unit disc. Probably, we are giving new information only for
what concerns   the case (2) but, for the sake of completeness, we
 deal with all  cases. It is interesting to
compare our approach to this embedding problem with the one given in
\cite[Section 5.9]{Shoiket}.

Given a point $a\in \mathbb{D}\setminus\{0\}$  there exists a unique
$\lambda \in \mathbb{C}$ such that $\func{Re}\lambda <0,$
$\func{Im}\lambda \in (-\pi ,\pi ]$ and $e^{\lambda }=a$. The {\sl
canonical spiral associated to $a$} is the curve
\begin{equation*}
\gamma_a:[1,\infty)\longrightarrow\D,\quad \gamma
_{a}(t):=e^{\lambda t}.
\end{equation*}%
Note that $\gamma _{a}$ is a spiral (actually a segment if $a\in
(0,1)$) which goes from $a$ to zero. The curve $\gamma _{a}$ has
finite length $\ell (\gamma _{a})$ given by
\begin{equation*}
\ell (\gamma _{a})=\int_{1}^{+\infty }|\gamma _{a}^{^{\prime }}(t)|dt=\frac{%
|\lambda a|}{\func{Re}(-\lambda )}\geq |a|.
\end{equation*}

\begin{proposition}\label{embedd}
Let $\varphi $ be an arbitrary element of $\mathrm{LFM}(\mathbb{D},\mathbb{D}%
).$

\begin{enumerate}
\item If $\varphi $ is trivial, neutral-elliptic, hyperbolic or
parabolic, then $\varphi $ can be always embedded into a semigroup
in$\mathbb{\ D}.$

\item If $\varphi $ is   attractive-elliptic  with
Denjoy-Wolff
point $\tau \in \mathbb{D}$ and repulsive fixed point $\beta \in \mathbb{C}%
_{\infty }\setminus \mathbb{D}$, let $\ell$ be the length of the
canonical spiral associated to $\varphi ^{\prime }(\tau )\in
\mathbb{D\setminus \{}0\mathbb{\}}.$ Then, $\varphi $ can be
embedded into a semigroup in$\mathbb{\ D}$ if and only if
\begin{equation}\label{condemb}
\left\vert \overline{\tau }-\dfrac{1}{\beta }\right\vert \ell\leq \left\vert \varphi ^{\prime }(\tau )\right\vert \left\vert 1-%
\dfrac{\tau }{\beta }\right\vert .
\end{equation}
[Again, $\infty ^{-1}\,$\ means $0$].
\end{enumerate}
\end{proposition}

\begin{proof}
$(1)$ If $\varphi $ is trivial or neutral-elliptic, then
$\v\in\Aut(\D)$ and, up to conjugation with an automorphism of $\D$
which maps the Denjoy-Wolff point of $\v$ to $0$, by the Schwarz
lemma $\v(z)=e^{ia} z$ for some $a\in\R$ and clearly $\v$ can be
embedded into the semigroup $(z,t)\mapsto e^{ita}z$. In the
hyperbolic or parabolic case, we can conjugate $\v$ with a Cayley
transform from $\D$ to $\H$ which maps the Denjoy-Wolff point of
$\v$ to $\infty$. Thus we obtain a linear fractional self-map $\phi$
of $\H$ of the form $w\mapsto e^aw+b$ with $a \in (0,+\infty )$ and
$b\in \mathbb{C}$,   $\func{Re}b\geq 0$ in the hyperbolic case; and
$a=0$, $b\in \mathbb{C}$ with $\func{Re}b\geq 0$ in the parabolic
case. Thus in the hyperbolic case $\phi$ belongs to the semigroup
$(w,t)\mapsto e^{a t}w+\frac{b}{e^a -1}(e^{a t}-1)$ while in the
parabolic case $\phi$ belongs to the semigroup $(w,t)\mapsto w+bt$.
Going back to the unit disc this implies that $\v$ can be embedded
into a semigroup in $\D$.

$(2)$ Let $\phi :=\alpha _{\tau }\circ \varphi \circ \alpha _{\tau
},$ with $\alpha _{\tau }(z)=\dfrac{\tau -z}{1-\overline{\tau
}z}$. The map $\phi $ is an attractive-elliptic linear fractional
map with   Denjoy-Wolff point  $0$ and  repulsive fixed point
$\widehat{\beta }=\alpha _{\tau }(\beta )\in \mathbb{C}_{\infty
}\setminus \mathbb{D}$ given by
\begin{equation*}
\phi (z)=\dfrac{az}{cz+1},
\end{equation*}%
where $a:=\varphi ^{\prime }(\tau )=\phi ^{\prime }(0)$, $0<|a|<1$ and $c=\dfrac{1%
}{\widehat{\beta }}(1-a)\in \mathbb{C}$. Let $\lambda\in\C$ be such
that $\func{Re}\lambda <0$, $\func{Im}\lambda \in (-\pi ,\pi ]$ and
$e^{\lambda }=a$. An easy computation shows that
\begin{equation}\label{condemb2}
\qquad \left\vert \overline{\tau }-\dfrac{1}{\beta }\right\vert \ell
\leq \left\vert \varphi ^{\prime }(\tau )\right\vert \left\vert
1-\dfrac{\tau }{\beta }\right\vert \Longleftrightarrow \left\vert \lambda \frac{c%
}{1-a}\right\vert \leq \func{Re}(-\lambda ).
\end{equation}
Now, assume \eqref{condemb} holds. By Theorem \ref{Tma Basico} and
\eqref{condemb2}, we deduce that
\begin{equation*}
F(z):=\lambda \frac{c}{1-a}z^{2}+\lambda z,\text{ }z\in \mathbb{D}
\end{equation*}
is the infinitesimal generator of a semigroup $(\phi _{t})$ of
$\mathrm{LFM}(\mathbb{D},\mathbb{D})$. Since $F(0)=0$ then $(\phi
_{t})$ is a non trivial-elliptic semigroup in $\mathbb{D}$. Using
Proposition \ref{Unival-VectorField}(1), one can find that the
Koenigs function of the semigroup is $h(z)=
\frac{z}{\frac{c}{1-a}z+1}$ and $\phi _{t}=h^{-1} (e^{\lambda
t}h).$ From this one can check that $\phi _{1}=\phi$. Thus  $\phi$
can be embedded in a semigroup in $\D$ and so does   $\varphi$.

Conversely, if $\v$ can be embedded in a semigroup in $\D$, so does
$\phi$. Let   $(\phi_t)$ be a semigroup in $\D$ such that
$\phi_1=\phi$. According to Theorem \ref{rigido}, $\phi_t$ is a
linear fractional self-map of $\D$ for all $t\geq 0$. If $h$ is the
Koenigs function of $(\phi_t)$ then $\phi _{t}(z)=h^{-1}(e^{\alpha
t}h(z))$ with $\func{Re}\alpha <0$ and
$h(z)=\dfrac{z}{\frac{c}{1-a}z+1}$ by \eqref{hK}. Using Proposition
\ref{Unival-VectorField}, one can deduce that the infinitesimal
generator of $(\phi _{t})$ is
\begin{equation*}
G(z)=\alpha \frac{c}{1-a}z^{2}+\alpha z,\text{ }z\in \mathbb{D}.
\end{equation*}%
By Theorem \ref{Tma Basico}
\begin{equation*}
\left\vert \alpha \frac{c}{1-a}\right\vert \leq \func{Re}(-\alpha
).
\end{equation*}%
Since $\phi _{1}=\phi $, $e^{\alpha }=\phi _{1}^{^{\prime }}(0)=\phi
^{\prime }(0)=a=e^{\lambda }.$ Therefore, $\lambda -\alpha =2k\pi
i,$ for some $k\in \mathbb{Z}.$ Bearing in mind that $\lambda \in
(-\pi ,\pi ],$ we find that $\func{Re}(\alpha )=\func{Re}(\lambda )$
and $|\alpha |\geq |\lambda |$ and by \eqref{condemb2} inequality
\eqref{condemb} holds.
\end{proof}

The above inequality explains dynamically well-known phenomena
concerning the embedding problem (shortly, EP). For instance, when
$\varphi ^{\prime }(\tau )\in (0,1),$ trivially the length of the
associated spiral
is exactly $\varphi ^{\prime }(\tau ).$ Since always $\left\vert \overline{%
\tau }-\dfrac{1}{\beta }\right\vert \leq \left\vert 1-\dfrac{\tau }{\beta }%
\right\vert ,$ we see that the answer to (EP) for those
attractive-elliptic maps is always positive. However, when $\varphi
^{\prime }(\tau )\in (-1,0),$ the length of the spiral is strictly
bigger than $\left\vert \varphi ^{\prime }(\tau )\right\vert ,$ so
we can have positive and negative answers. The inequality also says
that, for an attractive-elliptic element  of
$\mathrm{LFM}(\mathbb{D},\mathbb{D})$ the (EP) can be positively
solved whenever  the repulsive fixed point is close enough to
$\frac{1}{\overline{\tau }}$, namely when the map is similar enough
to a   neutral-elliptic map.

Moreover, if $\v$ is a hyperbolic or parabolic linear fractional
self-map of $\D$ and $\beta$ denote   the repelling fixed point of
$\v$ in the hyperbolic case and $\beta=\tau$ the Denjoy-Wolff point
of $\v$ in the parabolic case, it follows
\begin{equation*}
\left\vert \overline{\tau }-\dfrac{1}{\beta }\right\vert =\left\vert 1-%
\dfrac{\tau }{\beta }\right\vert \text{ and }\ell (\gamma
_{a})=\left\vert \varphi ^{\prime }(\tau )\right\vert =\varphi
^{\prime }(\tau )\in (0,1],
\end{equation*}
and thus \eqref{condemb} always holds for $\v$ hyperbolic or
parabolic.

Finally, if $\v$ is  neutral-elliptic,   since $\tau \in \mathbb{D}$
and the repulsive fixed point is $\frac{1}{\overline{\tau }}$ it
follows
\begin{equation*}
0=\left\vert \overline{\tau }-\overline{\tau }\right\vert \ell
(\gamma _{a})<\left\vert \varphi ^{\prime }(\tau )\right\vert
(1-|\tau |^{2}),
\end{equation*}
and even in this case \eqref{condemb} always holds.

\end{document}